\documentclass{article}
\usepackage{amsmath}[1996/11/01]
\usepackage{fullpage}
\usepackage{amssymb,amsthm,amsxtra}
\usepackage{epsfig}

\usepackage{color}
\definecolor{light-blue}{rgb}{0.8,0.85,1}
\definecolor{blue}{rgb}{0,0,1}
\definecolor{red}{rgb}{1,0,0}

\newtheorem{thm}{Theorem}
\newtheorem{cor}[thm]{Corollary}
\newtheorem{lem}[thm]{Lemma}
\newtheorem{prop}[thm]{Proposition}
\theoremstyle{remark}
\newtheorem*{rem}{Remark}

\newtheorem*{notrem}{\textbf{Notational Remark}}
\theoremstyle{definition}

\numberwithin{equation}{section}
\numberwithin{thm}{section}

\DeclareMathOperator{\Ai}{Ai}

\newcommand{\R}{\mathbb R}

\DeclareMathOperator{\erf}{erf}
\DeclareMathOperator{\Prob}{\mathbb P}

\begin{document}

\title{{\bf Painlev\'e formulas of the limiting distributions
for non-null complex sample covariance matrices}}
\author{{\bf Jinho Baik}
\footnote{Department of Mathematics, University of Michigan, Ann
Arbor, MI, 48109, baik@umich.edu}}

\date{July 20, 2005}
\maketitle

\begin{abstract}
In a recent study of large non-null sample covariance matrices, a
new sequence of functions generalizing the GUE Tracy-Widom
distribution of random matrix theory was obtained. This paper
derives Painlev\'e formulas of these functions and use them to
prove that they are indeed distribution functions.
Applications of these new distribution functions to last passage
percolation, queues in tandem and totally asymmetric simple
exclusion process are also discussed. As a part of the proof, a
representation of orthogonal polynomials on the unit circle in
terms of an operator on a discrete set is presented.
\end{abstract}


\section{Introduction}

Let $\Ai(u)$ denote the Airy function. It has an integral
representation
\begin{equation}\label{eq:Airy}
  \Ai(u) = \frac1{2\pi} \int e^{i(ua+\frac13 a^3)} da
\end{equation}
where the integral is over a curve from $\infty e^{5i\pi/6}$ to
$\infty e^{i\pi/6}$. The \emph{Airy kernel} (see, e.g.
\cite{Forrester, TW1}) is defined as
\begin{equation}\label{eq:Airyoperator}
  \mathbf{A}(u,v) := \frac{\Ai(u)\Ai'(v)-\Ai(u)\Ai'(v)}{u-v}
= \int_0^\infty \Ai(u+z)\Ai(z+v) dz.
\end{equation}
Let $\mathbf{A}_x$ be the \emph{Airy operator} acting on
$L^2((x,\infty))$ whose kernel is given by $\mathbf{A}(u,v)$.
Define
\begin{equation}\label{eq:F0op}
  F_0(x):= \det\bigl(1-\mathbf{A}_x\bigr).
\end{equation}
The `GUE Tracy-Widom distribution function' $F_0(x)$ is the
limiting distribution function of various models in mathematical
physics, probability and statistics (see e.g. \cite{TWICM} and
references in it).\footnote{In many literatures, $F_0$ is denoted
by $F_2$. In this paper, we reserve $F_2$ for a different
function.} Especially in statistics, the largest eigenvalue of the
sample covariance matrix of complex Gaussian samples with the
identity covariance (the so-called null case) is known to have the
limiting distribution given by $F_0(x)$. An intriguing result by
Tracy and Widom \cite{TW1} is that the Fredholm determinant has an
alternative expression:
\begin{equation}\label{eq:detP}
  F_0(x)= \det\bigl( 1- \mathbf{A}_x \bigr)
  = \exp \biggl( -\int_x^\infty (s-x)u^2(s) ds \biggr),
\end{equation}
where $u(x)$ is the solution to the \emph{Painlev\'e II} equation
\begin{equation}\label{eq:PII}
  u''=2u^3+xu,
\end{equation}
subject to the condition
\begin{equation}\label{eq:PIIboundary}
  u(x) \sim -\Ai(x) \qquad \text{as $x\to+\infty$}.
\end{equation}
It is known \cite{HM} that there is a unique global solution to
the equation \eqref{eq:PII} with the condition
\eqref{eq:PIIboundary}, and the solution satisfies (see, e.g.
\cite{HM, DZ2})
\begin{eqnarray}
  u(x) &=& -\Ai(x) + O\biggl( \frac{e^{-\frac14x^{3/2}}}{x^{1/4}} \biggr),
\qquad x\to +\infty \\
  u(x) &=& -\sqrt{\frac{-x}{2}} \biggl( 1+ O\biggl( \frac1{x^2} \biggr) \biggr),
\qquad x\to -\infty.
\end{eqnarray}
Recall that $\Ai(x) \sim e^{-\frac23 x^{3/2}}/(2\sqrt\pi x^{1/4})$
as $x\to +\infty$. The right-hand-side of \eqref{eq:detP} provides
a practical formula to plot the graph of $F_0$ numerically.


For $m=1,2,3,\dots$ and for complex numbers $w_1,w_2, \dots$,
define
\begin{equation}\label{eq:defs}
  s^{(m)}(u; w_1, \dots, w_m) = s^{(m)}(w_1, \dots, w_m) :=
\frac1{2\pi} \int  e^{\frac13 ia^3+iua}
  \prod_{j=1}^m \frac1{w_j+ia} da
\end{equation}
where the contour is from $\infty e^{5i\pi/6}$ to $\infty
e^{i\pi/6}$ such that the poles $a=iw_1, \dots, iw_m$ lie
\emph{above} the contour. Also define
\begin{equation}\label{eq:tdef}
  t^{(m)}(v; w_1,\dots, w_{m-1}) =  t^{(m)}(w_1,\dots, w_{m-1})
:= \frac1{2\pi} \int  e^{\frac13 ib^3+ivb}
  \prod_{j=1}^{m-1} (w_j-ib) db
\end{equation}
where the contour is from $\infty e^{5i\pi/6}$ to $\infty
e^{i\pi/6}$. Comparing with~\eqref{eq:Airy}, $t^{(m)}$ is a sum of
derivatives of the Airy function. On the other hand, when
$w_1=\dots =w_m=0$, $s^{(m)}$ is a sum of anti-derivatives of the
Airy function. However for general $w_j$'s, $s^{(m)}$ is a
Cauchy-type transform of the integrand of the Airy function.
Define
\begin{equation}\label{eq:defFk}
\begin{split}
   & F_k(x; w_1, \dots, w_k) \\
&\quad := F_0(x) \cdot \det \biggl (\delta_{mn} -
<\frac1{1-\mathbf{A}_x} s^{(m)}(w_1, \dots, w_m), t^{(n)}(w_1,
\dots, w_{n-1})>_{L^2((x,\infty))} \biggr)_{1\le m,n\le k}
\end{split}
\end{equation}
where $<,>_{(x,\infty)}$ denotes the \emph{real} inner product in
$L^2((x,\infty))$;
\begin{equation}
\begin{split}
   &<\frac1{1-\mathbf{A}_x} s^{(m)}(w_1, \dots, w_m), t^{(n)}(w_1,
\dots, w_{n-1})>_{L^2((x,\infty))} \\
&\qquad = \int_x^\infty \biggl( \frac1{1-\mathbf{A}_x}
s^{(m)}(w_1, \dots, w_m)\biggr)(u) t^{(n)}(u; w_1, \dots, w_{n-1})
du.
\end{split}
\end{equation}
(It is well-known that $1-\mathbf{A}_x$ is invertible.) Set
\begin{equation}\label{eq:defFkk}
  F_k(x):= F_k(x; 0,0\dots, 0), \qquad k=1,2,\dots.
\end{equation}

The functions $F_k(x;w_1, \dots, w_k)$ were introduced recently in
\cite{BBP} as limits of the distribution functions of the largest
eigenvalues of certain non-null complex sample covariance matrices
and also other probability models. See Section~\ref{sec:models}
below for more details on the motivations. The purpose of this
paper is to find a Painlev\'e type formula for $F_k(x; w_1,
\cdots, w_k)$ analogous to \eqref{eq:detP}. Such formula is used
to prove that $F_k(x; w_1, \cdots, w_k)$ is indeed a distribution.
It also allows us to be able to plot the graph of $F_k(x)$.


\subsection{Results}

\subsubsection{Alternative determinantal formula}

We first obtain an alternative determinantal formula of $F_k(x;
w_1, \dots, w_k)$. The definition \eqref{eq:defFk} involves the
functions $s^{(m)}$ and $t^{(m)}$ and it is not transparent that
the formula is symmetric in $w_1, \dots, w_k$, which should be the
case from its origin in the sample covariance matrix \cite{BBP}
(see also Section \ref{sec:models} below). This symmetry is clear
in the following theorem.

For a complex number $w$, set
\begin{equation}\label{eq:Cdef}
  C_w(u) := \frac1{2\pi} \int  e^{i(\frac13 a^3+ua)}
  \frac1{w+ia} da
\end{equation}
where the contour is, as in the definition \eqref{eq:defs} of
$s^{(m)}$, from $\infty e^{5i\pi/6}$ to $\infty e^{i\pi/6}$ such
that the pole $a=iw$ lies above the contour. Hence
$s^{(1)}(u;w_1)=C_{w_1}(u)$. Also note that $t^{(1)}(v)=\Ai(v)$.

\begin{thm}\label{thm:det}
With above notations, for real $x$ and complex $w$ set
\begin{equation}\label{eq:fdef}
  f(x,w) := 1- < \frac1{1-\mathbf{A}_x} C_w, \Ai>_{L^2((x,\infty))}
  = 1- \int_x^\infty \biggl( \frac1{1-\mathbf{A}_x} C_w\biggr)(u) \Ai(u) du.
\end{equation}
For distinct complex numbers $w_1,\dots, w_k$,
\begin{equation}\label{eq:detformula}
\begin{split}
  F_k(x;w_1,\dots, w_k) = F_0(x) \cdot \frac{
  \det \begin{pmatrix}
  (w_m+ D_x)^{n-1}f(x,w_m)
  \end{pmatrix}_{1\le m,n\le k}}
  {\displaystyle \prod_{1\le m<n\le k} (w_n-w_m)}
\end{split}
\end{equation}
where $D_x=\frac{\partial}{\partial x}$ denotes the derivative
with respect to $x$. When some of $w_j$'s coincide, the above
formula still holds by using the l'Hosptial's rule for the
right-hand-side of \eqref{eq:detiden}.
\end{thm}

\begin{rem}
P. Deift and A. Its pointed out that this formula resembles the
Darboux transformation in the theory of integrable systems (see
e.g., \cite{Darboux}). It would be interesting to identify the
above formula in terms of a Darboux transformation of an
integrable system.
\end{rem}

This theorem follows from row and column operations
of~\eqref{eq:defFk} exploiting the fact that $t^{(n)}$ is a sum of
derivatives of $\Ai$ and that $s^{(m)}$ is a linear combination of
$C_{w_j}$. The proof is given in Section~\ref{sec:proofthmdet}.

\subsubsection{Painlev\'e formula}

In the next theorem, we show that the function $f(x,w)$ defined in
\eqref{eq:fdef} is related to the Painlev\'e II equation.

First we need a definition.
Let $M(z; x)= \biggl( \begin{smallmatrix} M_{11}(z) & M_{12}(z) \\
M_{21}(z) & M_{22}(z)
\end{smallmatrix} \biggr)$ be the $2\times 2$ matrix-valued solution to the
following Riemann-Hilbert problem:
\begin{itemize}
\item $M(z;x)$ is analytic for $z\in\mathbb{C}\setminus \mathbb{R}$
and is continuous for $z\in \overline{\mathbb{C}\setminus
\mathbb{R}}$
\item For $z\in \mathbb{R}$,
\begin{equation}\label{eq:jump}
   M_+(z;x) = M_-(z;x) \begin{pmatrix} 1 & -e^{-2i(\frac43
z^3+xz)} \\ e^{2i(\frac43 z^3+xz)} &0
\end{pmatrix}
\end{equation}
where $M_+(z;x)$ (resp. $M_-(z;x)$) denotes the limit of $M(z';x)$
as $z'\to z$ from the bottom (resp. top) of the contour
$\mathbb{R}$.
\item $M(z;x) \to I$ as $z\to\infty$.
\end{itemize}
The precise statement of the last condition is the following:
there is $\epsilon>0$ such that $M(z;x)= I+O(z^{-1})$ uniformly as
$z\to\infty$ for $z$ in sectors $\epsilon<Arg(z)<\pi -\epsilon$
and $\pi +\epsilon<Arg(z)<2\pi -\epsilon$, and $M(z;x)$ is bounded
for all $z\in \mathbb{C}\setminus \mathbb{R}$.

This is the Riemann-Hilbert problem for the Painlev\'e II equation
when the so-called monodromy data satisfies $p=-q=1$ and $r=0$
\cite{IN, FZ, DZ2}. It is known that there is a unique solution to
this Riemann-Hilbert problem. Moreover, as $z\to\infty$, there is
an expansion of form
\begin{equation}
  M(z;x) = I + \frac{M_1(x)}{z} + O \biggl( \frac1{z^2} \biggr),
\qquad M_1(x) = \frac1{2i} \begin{pmatrix}
  -v(x) & u(x) \\ -u(x) & v(x)
\end{pmatrix}
\end{equation}
where $u(x)$ is the solution of the Painlev\'e II equation
\eqref{eq:PII} satisfying \eqref{eq:PIIboundary}, and
\begin{equation}
  v(x) = \int_\infty^x u(s)^2 ds.
\end{equation}

The following theorem shows that $f(x,w)$ is expressible in terms
of the Riemann-Hilbert problem for Painlev\'e II equation.
\begin{thm}\label{thm:fPain}
The function $f(x,w)$ defined in \eqref{eq:fdef} satisfies the
following:
\begin{equation}\label{eq:fPain}
  f(x,w) =
  \begin{cases} M_{22}(-\frac12iw;x), \qquad & Re(w)>0 \\
  -M_{21}(-\frac12iw;x)e^{\frac13 w^3-xw}, \qquad & Re(w)<0.
  \end{cases}
\end{equation}
\end{thm}

Note that from the jump condition \eqref{eq:jump}, $f(x,w)$ is
continuous for $w\in\mathbb{R}$, and hence is an entire function
in $w$.

Together with Theorem~\ref{thm:det}, Theorem~\ref{thm:fPain}
yields the desired Painlev\'e II formula of $F_k(x;w_1,\dots,
w_k)$, which is the main result of this paper.

\begin{cor}\label{cor:main}
The function $F_k(x;w_1,\dots, w_k)$ defined by~\eqref{eq:defFk}
is equal to~\eqref{eq:detformula} with $f(x,w)$ given
by~\eqref{eq:fPain}.
\end{cor}

The function given in the right-hand-side of~\eqref{eq:fPain} had
previously appeared in \cite{BR2} (equation (2.22)) and \cite{BR4}
(equation (3.5)) as a limiting function for a last passage site
percolation model. In the context of symmetrized random
permutations and last passage percolation models, \cite{BR2, BR4}
showed, among other things, the $k=1$ case of
Corollary~\ref{cor:main};
\begin{equation}
    F_1(x, w_1)= F_0(x) f(x,w_1)
\end{equation}
where $f(x,w)$ given by the right-hand-side of~\eqref{eq:fPain}.
This paper proves that the general case is expressible in terms of
derivatives of the same function $f(x,w)$.

\subsubsection{Properties of $f(x.w)$}

The papers \cite{BR2, BR4} proved several properties of the
function defined by the right-hand-side of~\eqref{eq:fPain}. By
setting $w\mapsto \frac12w$ in Lemma 2.1 of \cite{BR2} or Lemma
3.1 of \cite{BR4}, we find the following properties of $f(x,w)$.
The following complementary function is useful: set
\begin{equation}
   g(x,w) :=
  \begin{cases} M_{12}(-\frac12iw;x), \qquad & Re(w)>0 \\
  -M_{11}(-\frac12iw;x)e^{\frac13 w^3-xw}, \qquad & Re(w)<0.
  \end{cases}
\end{equation}

\begin{lem}[\cite{BR2, BR4}]\label{lem:Pain} The following holds.
\begin{enumerate}
\item $f(x,w), g(x,w)$ are real for $w\in \mathbb{R}$.
\item For each fixed $w\in\mathbb{C}$, as $x\to +\infty$
\begin{eqnarray}
  f(x,w)&=& 1+O(e^{-cx^{3/2}}),\\
  g(x,w)&=& -e^{\frac13w^3-xw} \bigl( 1+O(e^{-cx^{3/2}})\bigr)
\end{eqnarray}
and as $x\to -\infty$,
\begin{eqnarray}\label{eq:fneginfty}
  f(x,w) &\sim&
\frac1{\sqrt{2}} e^{\frac16 w^3-\frac{1}{6}|x|^{3/2}+\frac12|x|w
-w^2|x|^{1/2}} \\
  g(x,w) &\sim&
-\frac1{\sqrt{2}} e^{\frac16 w^3-\frac1{6}|x|^{3/2}+\frac12|x|w
-w^2|x|^{1/2}}.
\end{eqnarray}
\item
\begin{eqnarray}
\lim_{w\to+\infty} f(x,w)=1, && \lim_{w\to+\infty} g(x,w)=0,\\
\lim_{w\to-\infty} f(x,w)=0, && \lim_{w\to-\infty} g(x,w)=0,\\
\label{eq:initial} f(x,0)= \mathcal{E}(x),  &&
g(x,0)=-\mathcal{E}(x),
\end{eqnarray}
where
\begin{equation}
  \mathcal{E}(x) := \exp\biggl\{ \int_x^\infty u(s)ds \biggr\}.
\end{equation}
($\mathcal{E}(x)$ was denoted by $E(x)^2$ in \cite{BR2, BR4}.)
\item For all $x\in \mathbb{R}$ and $w\in\mathbb{C}$,
\begin{eqnarray}\label{eq:Lax1}
  \frac{\partial}{\partial x} \binom{f(x,w)}{g(x,w)} &=& \begin{pmatrix} 0 & u(x) \\
  u(x) & -w \end{pmatrix} \binom{f(x,w)}{g(x,w))},\\
\label{eq:Lax2}  \frac{\partial}{\partial w}
\binom{f(x,w)}{g(x,w)} &=&
\begin{pmatrix} (u(x))^2 & -wu(x)-u'(x) \\
-wu(x)+u'(x) & w^2-x-(u(x))^2 \end{pmatrix}
\binom{f(x,w)}{g(x,w)}.
\end{eqnarray}
\item
\begin{eqnarray}
  f(x,w)&=&-g(x,-w) e^{\frac13w^3-xw}, \\
  g(x,w)&=&-f(x,-w) e^{\frac13w^3-xw}.
\end{eqnarray}
\item For each fixed $y\in\R$, as $w\to-\infty$,
\begin{eqnarray}
  f(y\sqrt{|w|}+w^2,w) &\to& \erf(y)=\frac1{\sqrt{2\pi}} \int_{-\infty}^y e^{-\frac12 s^2}ds,\\
g(y\sqrt{|w|}+w^2,w) &\sim&
-e^{\frac{2}3|w|^3+\sqrt{2}y|w|^{3/2}}.
\end{eqnarray}
\end{enumerate}
\end{lem}

Note that \eqref{eq:Lax1} and \eqref{eq:Lax2} are the Lax pair
equations for the Painlev\'e II equation. Hence
Theorem~\ref{thm:det}, Theorem~\ref{thm:fPain} and
Lemma~\ref{lem:Pain} yields that $F_k(x;w_1,\dots, w_k)$ is
expressible in terms of the Lax pair equations of the Painlev\'e
II equation.

\begin{rem}
After this paper was completed, Harold Widom found a different
proof of \eqref{eq:Lax1} and \eqref{eq:Lax2} for $f(x,w)$ defined
by \eqref{eq:fdef} and $g(x,w)$ defined by $g(x,w) = -\bigl(
\frac1{1-\mathbf{A}_x} C_w\bigr)(x)$ using the method of
\cite{TW1}. The proof of Widom is algebraic and is more direct. On
the other hand, the current paper proves a general identity and
then takes a limit as outlined in subsection \ref{sec:outline}
below.
\end{rem}

From \eqref{eq:Lax1} and \eqref{eq:Lax2}, $f(x,w)$ itself
satisfies a second order linear differential equations in $x$ and
$w$ with coefficients involving $u(x)$.

\begin{cor}
Denoting $\frac{\partial}{\partial x} f(x,w) = f'(x,w)$ and
$\frac{\partial}{\partial w} f(x,w) = \dot{f}(x,w)$, $f$ satisfies
\begin{equation}
  -f''  + \biggl( \frac{u'}{u} -w \biggr) f'  +u^2 f=0
\end{equation}
and
\begin{equation}\label{eq:DE}
  - \ddot{f} + \biggl( \frac{u}{wu+u'} + w^2-x \biggr) \dot{f} +
  \biggl( - \frac{u^3}{wu+u'} + u^4+xu^2-(u')^2 \biggr) f=0.
\end{equation}
\end{cor}
Together with the initial conditions $f(x,0)=\mathcal{E}(x)$ and
$\dot{f}(x,0)=(u^2(x)+u'(x))\mathcal{E}(x)$, \eqref{eq:DE}  may
provide a numerical way to compute the function $f(x,w)$, and
hence $F_k(x;w_1,\dots, w_k)$.

\subsubsection{Formula of $F_k(x)$}

\begin{figure}
 \centerline{\epsfig{file=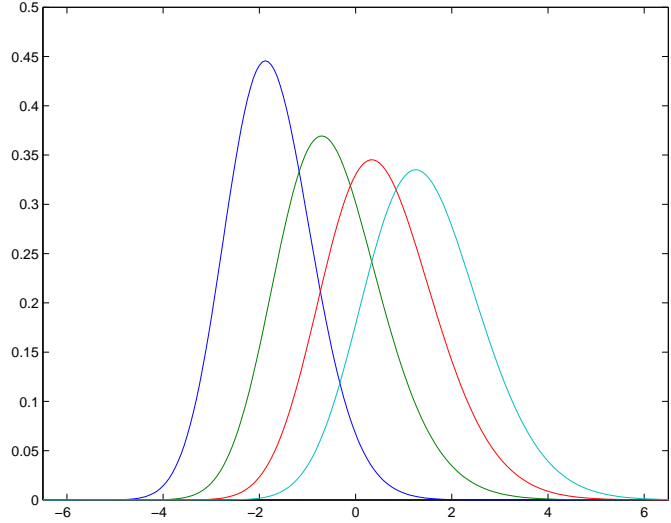, height=7cm}}
 \caption{Graph of $\frac{d}{dx}F_j(x)$, $j=0,1,2,3$ from the left to the right}\label{fig:inter}
\end{figure}

When $w_1=\dots =w_k=0$, using the l'Hospitals' rule
in~\eqref{eq:detformula},
\begin{equation}\label{eq:detFwhenzero}
  F_k(x) = F_k(x; 0,0,\dots, 0)= \frac1{\prod_{j=0}^{k-1} j!} F_0(x) \cdot
  \det \biggl( D_w^{m-1} \bigl\{ (w+D_x)^{n-1}
  f(x,w) \bigr\} \bigr|_{w=0} \biggr)_{1\le m,n\le k} .
\end{equation}
By using \eqref{eq:Lax1}, \eqref{eq:Lax2} and \eqref{eq:initial},
one can in principle compute the determinant. The first three
functions are
\begin{equation}\label{eq:FPain}
\begin{split}
  F_1(x) &= F_0(x) \mathcal{E}(x),\\
  F_2(x)  &= F_0(x) \mathcal{E}(x)^2 \bigl\{ 1+ u(x+2u^2+2u') \bigr\},\\
  F_3(x) &= F_0(x) \mathcal{E}(x)^3\bigl\{ 1+ 2u(x+2u^2+2u') + \frac12(u^2-u') (x+2u^2+2u')^2
  \bigr\}.
\end{split}
\end{equation}
Using the numerical evaluation of the Painlev\'e solution $u(x)$
which is available at the website of M. Pr\"ahofer
(\verb+http://www-m5.ma.tum.de/KPZ+), these formulas provide a
convenient way to plot the graphs of $F_k$. Figure~\ref{fig:inter}
is the graphs of the density function $\frac{d}{dx}F_j(x)$ for
$j=0,1,2,3$.
Note that the function moves to the right as the index $k$
increases. The numerical means and the standard deviations of
$F_k(x)$ are the following:

\vspace{0.3cm}

{\centerline{
\begin{tabular}{|c|c|c|}
  \hline
   & mean  & standard deviation \\ \hline
  $F_0$ & -1.771\dots & 0.90\dots \\ \hline
  $F_1$ & -0.494\dots & 1.11\dots \\ \hline
  $F_2$ & 0.543\dots & 1.18\dots \\ \hline
  $F_3$ & 1.445\dots & 1.21\dots \\ \hline
\end{tabular}
}}

\subsubsection{$F_k$ are distribution functions}

The Painlev'e formula obtained above can be used to prove that
$F_k(x;w_1, \dots, w_k)$ is indeed a distribution function.

\begin{cor}
The function $F_k(x;w_1, \dots, w_k)$ is a distribution for real
$w_1, \dots, w_k$.
\end{cor}

\begin{proof}
In \cite{BBP}, the function $F_k(x;w_1, \dots, w_k)$ are shown to
be continuous, non-decreasing and converges to $1$ as $x\to
+\infty$ (see the paragraph after (25)). We need to show that
$F_k(x;w_1, \dots, w_k)\to 0$ as $x\to -\infty$. From
\eqref{eq:fneginfty}, all entries of the determinant in both
\eqref{eq:detformula} and \eqref{eq:detFwhenzero} are in absolute
value less than or equal to $Ce^{-c|x|^{3/2}}$ for some constants
$C, c>0$ as $x\to -\infty$. Also $F_0(x)\le Ce^{-c|x|^3}$ for some
other constants $C, c>0$ (see e.g. (2.13) of \cite{BR2}). Hence
$F_k(x; w_1, \dots, w_k) = O(e^{-c|x|^3})$ for some constant $c>0$
as $x\to -\infty$.
\end{proof}

\subsection{Outline of the proof and orthogonal polynomials on the
unit circle}\label{sec:outline}

Theorem \ref{thm:det} is obtained by applying a sequence of row
and column operations to the determinant \eqref{eq:defFk}. This
part is the main bulk of the paper and the proof is given in
Section \ref{sec:proofthmdet}.

The proof of Theorem \ref{thm:fPain} is indirect. We use a
representation of orthogonal polynomials on the unit circle in
terms of an operator on a discrete set. Since such a
representation may be interesting in itself, we state it here.
This formula follows from a general identity (see~\eqref{eq:GCBO}
below) between Toeplitz determinants and Fredholm determinants on
integer lattice obtained by Geronimo and Case \cite{GeronimoC69},
and also independently by Borodin and Okounkov \cite{BoOk} (see
also \cite{BW, BottSimple} for shorter proofs).

Set $\mathbb{N}_0:=\mathbb{N}\cup\{0\}$. Let $\phi(z)$ be a
function which is positive on the unit circle. For simplicity of
argument, we assume that $\phi(z)$ is analytic in neighborhood of
the unit circle. 
Let $\phi=\phi_+\phi_-$ be a Wiener-Hopf factorization of $\phi$
where $\phi_+$ extends to a non-vanishing analytic function
interior of the circle and $\phi_-$ extends to a non-vanishing
analytic function exterior of the circle.
Set
\begin{equation}
  \psi(z):= \frac{\phi_+(z)}{\phi_-(z)}.
\end{equation}
For a function $f(z)$ on the unit circle, $f_k$
denotes its $k^{th}$ Fourier coefficient:
\begin{equation}
  f_k := \int_{|z|=1} z^{-k} f(z) \frac{dz}{2\pi iz}.
\end{equation}

\begin{prop}\label{prop:pidet}
Let $\pi_n(z)$ be the monic orthogonal polynomial on the unit
circle with respect to the measure $\phi(z) \frac{dz}{2\pi iz}$,
and let $\pi_n^*(z)=z^n\overline{\pi}_n(\frac1{z})$ be its
$*$-transform. For $\phi$ satisfying above conditions,
\begin{equation}\label{eq:pidet}
  \pi_n^*(z) = 
   e^{-\sum_{k=1}^\infty (\log \phi)_{k} z^k}
   \biggl\{ 1- <
   \frac1{1-P_n A B P_n} P_nQ, P_nR>_{\ell^2(\mathbb{N}_0)}
   \biggr\}, \qquad |z|<1,
\end{equation}
where $<,>_{\ell^2(\mathbb{N}_0)}$ is the real inner product on
$\ell^2(\mathbb{N}_0)$, $P_n$ is the projection on the set $\{n,
n+1, n+2, \dots\}$, the operators $A, B : \ell^2(\mathbb{N}_0) \to
\ell^2(\mathbb{N}_0)$ are defined by the kernels
\begin{equation}\label{eq:defAB}
  A(j,m) = (\psi^{-1})_{j+m+1}, \qquad
  B(m,k) = \psi_{-m-k-1}
\end{equation}
and the functions $Q, R\in \ell^2(\mathbb{N}_0)$ are given by
\begin{equation}\label{eq:defQR}
  Q(j) = (\psi^{-1})_{j+1}, 
  \qquad
  R(k) =  \biggl(
  \frac{z}{\cdot-z} \psi(\cdot) \biggr)_{-k-1}
  = \int_{|b|=1} b^{k+1} 
  \frac{z}{b-z} \psi(b)\frac{db}{2\pi ib}.
\end{equation}
On the other hand,
\begin{equation}\label{eq:pidet222}
  \pi_n(z) = 
   z^ne^{-\sum_{k=1}^\infty (\log \phi)_{-k} z^{-k}}
   \biggl\{ 1- <
   \frac1{1-P_n A B P_n} P_nU, P_nV>_{\ell^2(\mathbb{N}_0)}
   \biggr\}, \qquad |z|>1,
\end{equation}
where
\begin{equation}\label{eq:defUV}
  U(j) =  \biggl(
  \frac{\cdot}{z-\cdot} \psi^{-1}(\cdot) \biggr)_{j+1}
  = \int_{|a|=1} a^{-j-1} 
  \frac{a}{z-a} \frac1{\psi(a)}\frac{da}{2\pi ia},
  \qquad V(k)=\psi_{-k-1},
\end{equation}
\end{prop}

\begin{rem}
(i) As $\pi_n^*(z)$ is an entire function, the formula
\eqref{eq:pidet} also holds for a region of $|z|\ge 1$ to which
the right-hand-side of \eqref{eq:pidet} is analytically continued.
(ii) The conditions for $\phi$ above can be weakened, but we do
not discuss such an issue in this paper.
(iii) Recall that Wiener-Hopf factorizations of $\phi$ are
different by a factor of a multiplicative constant. However, since
$A$ and $B$ have both factors $\psi^{-1}$ and $\psi$,
respectively, the operator $AB$ is unaffected by a different
choice of $\phi_+$ and $\phi_-$. The inner products in
\eqref{eq:pidet} and \eqref{eq:pidet222} also remain the same even
if $\psi$ is multiplied by a constant. Therefore, \eqref{eq:pidet}
and \eqref{eq:pidet222} does not depend on the choice of a
Wiener-Hopf factorization of $\phi$.
\end{rem}

\begin{rem}
After this paper was completed, it turned out during a
conversation with Andrei Mart\'{i}nez-Finkelstein that
\eqref{eq:pidet222} appeared in \cite{MMS} in a very different
form. In the first formula of the equation (32) of \cite{MMS}, the
authors found a series expansion for $\pi(z)$. However one can
check that the series is precisely what one would obtain once the
Neuman series of the operator $P_nABP_n$ is taken in
\eqref{eq:pidet222}. The proof of \cite{MMS} is based on a
Riemann-Hilbert method. By turning the argument backward, it is
also possible to prove the identity of Geromino-Case and
Borodin-Okounkov using a Riemann-Hilbert method. This will appear
in a future work.
\end{rem}

We regard \eqref{eq:pidet} as an identity. We take a special
choice of $\phi$ and then take a limit of both sides of the
identity \eqref{eq:pidet}. A steepest-descent analysis shows that
the right-hand-side converges to the formula \eqref{eq:fdef}. On
the other hand, a Riemann-Hilbert asymptotic analysis to the
left-hand-side yields the Painlev\'e formula \eqref{eq:fPain}.
Hence the identity \eqref{eq:fPain} follows from the identity
between the orthogonal polynomials and their operator
representation.


%

\medskip

This paper is organized as follows. In Section~\ref{sec:models},
we present the statistical and probabilistic models in which the
distributions $F_k(x;w_1,\dots, w_k)$ appear. The proof of Theorem
\ref{thm:det} is given in Section~\ref{sec:proofthmdet}.
Section~\ref{sec:ProofPropdidet} proves
Proposition~\ref{prop:pidet} and
Section~\ref{sec:proofTheoremfPain} proves
Theorem~\ref{thm:fPain}.

\medskip
\noindent {\bf Acknowledgments.} The authors would like to thank
P. Deift, A. Its, A. Mart\'{i}nez-Finkelshtein, B. Simon and H.
Widom for useful conversations and communications. This work was
supported in part by NSF Grant \#DMS-0350729 and the AMS
Centennial Fellowship.

\section{Models}\label{sec:models}

We discuss several statistics and probability models in which
$F_k$'s appear.

\subsection{Non-null complex sample covariance
matrices}\label{sec:samplecovariance}

Let $M\ge N\ge 1$ be integers. Let $\vec{y}_1, \dots, \vec{y}_M$
be independent \emph{complex} Gaussian $N\times 1$ column vectors
with mean $\vec{\mu}$ and population covariance $\Sigma$: the
density of $\vec{y}_1$ is
\begin{equation}
  p(\vec{y}_1) = \frac1{(2\pi)^{N/2}(\det\Sigma)^{1/2}}
  e^{-\frac12(\vec{y}_1-\vec{\mu})^*\Sigma^{-1}
  (\vec{y}_1-\vec{\mu})}
\end{equation}
where $*$ denotes the complex transpose. Denote by $\bar{Y}$ the
sample mean $\bar{Y}:= \frac{1}{M} (\vec{y}_1+\dots+\vec{y}_M)$
and by $X=[\vec{y}_1-\bar{Y}, \dots, \vec{y}_M-\bar{Y}]$ the
(centered) $M\times N$ sample matrix. Define the \emph{sample
covariance matrix} by
\begin{equation}
  S = \frac1{M} XX^*.
\end{equation}
When the covariance matrix $\Sigma$ is the identity matrix, the
distribution of the eigenvalues of $S$ is sometimes called the
Laguerre unitary ensemble and is well-studied in the random matrix
theory (see e.g. \cite{ForresterBook}). In particular, as
$M,N\to\infty$ while $M/N=\gamma^2$ is in a compact subset of
$[1,\infty)$, the largest eigenvalue $\lambda_{\max}$ satisfies
the limit law (see e.g. \cite{Forrester, kurtj:shape})
\begin{equation}\label{eq:F0limitlaw}
  \Prob\biggl( \bigl(\lambda_{\max}- (1+\gamma^{-1})^2)\bigr)
  \cdot \frac{\gamma}{(1+\gamma)^{4/3}}M^{2/3} \le x \biggr)
  \to F_0(x),
\end{equation}
where $F_0(x)$ is the Tracy-Widom distribution~\eqref{eq:F0op}.

Johnstone \cite{Johnstone} proposed the study of the so-called the
`spiked population model' where the covariance matrix $\Sigma$ is
a finite rank perturbation of the identity matrix. For possible
applications of the spiked population model in statistics, finance
and telecommunications, see the references in \cite{Johnstone} and
\cite{BBP}. For spiked population models, it is interesting to
determine the effect of non-unit eigenvalues of the covariance
matrix on the largest eigenvalue of the sample covariance matrix.
For complex Gaussian samples, \cite{BBP} determined the critical
value of the non-unit covariance eigenvalue. When some of the
non-unit eigenvalues of the covariance matrix are above the
critical value, $\lambda_{\max}$ behaves differently
from~\eqref{eq:F0limitlaw}. The function $F_k(x)$ is the limiting
distribution of the $\lambda_{\max}$ when the largest eigenvalue
of the covariance matrix is of multiplicity $k$ and is equal to
the critical value.

 Let $\ell_1\ge \dots\ge \ell_r>0$ be the
non-unit eigenvalues of $\Sigma$ where $r$ is independent of $M$
and $N$.

\begin{thm}[Theorem 1.1 of \cite{BBP}]\label{thm:BBP1}
As $M, N\to\infty$ such that $M/N=\gamma^2$ lies in a compact
subset of $[1,\infty)$, the following holds.
\begin{itemize}
\item[(a)] When
\begin{equation}
  \ell_1=\dots =\ell_k= 1+\gamma^{-1}
\end{equation}
for some $0\le k\le r$, and $\ell_{k+1}, \dots, \ell_r$ are in a
compact subset of $(0, 1+\gamma^{-1})$,
\begin{equation}
  \Prob\biggl( \bigl(\lambda_{\max}- (1+\gamma^{-1})^2)\bigr)
  \cdot \frac{\gamma}{(1+\gamma)^{4/3}}M^{2/3} \le x \biggr)
  \to F_k(x)
\end{equation}
where $F_k(x)$ is defined in \eqref{eq:defFkk}.
\item[(b)] When
\begin{equation}
  \text{$\ell_1=\dots = \ell_k$ are in a compact subset of
  $(1+\gamma^{-1}, \infty)$}
\end{equation}
for some $1\le k\le r$, and $\ell_{k+1}, \dots, \ell_r$ are in a
compact subset of $(0, \ell_1)$,
\begin{equation}
  \Prob\biggl( \bigl(\lambda_{\max}- (\ell_1+\frac{\ell_1\gamma^{-2}}{\ell_1-1})
  \bigr)
  \cdot \sqrt{M} \sqrt{\ell_1^2-\frac{\ell_1^2\gamma^{-2}}{(\ell_1-1)^2}}
  \le x \biggr)
  \to G_k(x)
\end{equation}
where $G_k(x)$ is the distribution of the largest eigenvalue of
$k\times k$ Gaussian unitary ensemble.
\end{itemize}
\end{thm}

More detailed nature of the phase transition around the critical
value $1+\gamma^{-1}$ was also studied in the same paper.

\begin{thm}[Theorem 1.2 of \cite{BBP}]\label{thm:BBP2}
For some $1\le k\le r$, set
\begin{equation}
  \ell_j= 1+\gamma^{-1} - \frac{(1+\gamma)^{3/2}w_j}{\gamma
  M^{1/3}},
  \qquad j=1,2,\dots, k.
\end{equation}
When $w_1, \dots, w_k$ are in a compact subset of $\mathbb{R}$ and
$\ell_{k+1}, \dots, \ell_r$ are in a compact subset of $(0,
1+\gamma^{-1})$, as $M, N\to\infty$ while $M/N=\gamma^2$ is in a
compact subset of $[1, \infty)$,
\begin{equation}
  \Prob\biggl( \bigl(\lambda_{\max}- (1+\gamma^{-1})^2)\bigr)
  \cdot \frac{\gamma}{(1+\gamma)^{4/3}}M^{2/3} \le x \biggr)
  \to F_k(x; w_1, \dots, w_k)
\end{equation}
where $F_k(x; w_1, \dots, w_k)$ is defined in \eqref{eq:defFk}.
\end{thm}

It is transparent from this theorem that $F_k(x;w_1, \dots, w_k)$
should be symmetric in $w_1, \dots, w_k$ since re-labelling the
eigenvalues does not change the limit law. Further work on the
eigenvalues of the spiked model can be found in \cite{DPaul,
BSilverstein}.

\subsection{Last passage percolation and queues in tandem}

Suppose that to each lattice points $(i,j)\in \mathbb{Z}^2$, an
independent random variable $X(i,j)$ is associated. Let
$(1,1)\nearrow (N, M)$ denote the set of `up/right paths'
$\pi=\{(i_k, j_k)\}_{k=1}^{N+M-1}$ where $(i_{k+1},
j_{k+1})-(i_k,j_k)$ is either $(1,0)$ or $(0,1)$, and
$(i_1,j_1)=(1,1)$ and $(i_{N+M-1}, j_{N+M-1})=(N, M)$. Note that
the cardinality of $(1,1)\nearrow (N,M)$ is $\binom{N+M-2}{N-1}$.
Set
\begin{equation}
  L(N, M) := \max_{\pi \in (1,1)\nearrow (N,M)} \sum_{(i,j)\in
  \pi} X(i,j).
\end{equation}
By interpreting $X(i,j)$ as the (random) time spent to pass the
site $(i,j)$, $L(N,M)$ is the \emph{last passage time} to travel
from $(1,1)$ to $(N,M)$ along an admissible up/right path.

Recall that the exponential random variable of mean $m$ has the
density function $\frac1{m} e^{-x/m}$, $x\ge 0$. It is known that
(see e.g. Proposition 6.1 of \cite{BBP}; we here scale $X(i,j)$ of
\cite{BBP} by $M$) when $X(i,j)$ is an exponential random variable
of mean $\ell_i$ (independent of $j$), $\frac{L(N,M)}{M}$ has the
same distribution as the largest sample eigenvalue
$\lambda_{\max}$ of complex Gaussian samples when the eigenvalues
of the population covariance matrix $\Sigma$ are $\ell_1, \dots,
\ell_N$. Therefore for the last passage percolation model which
have the identically distributed passage time for all but finitely
many columns, Theorems~\ref{thm:BBP1} and Theorem~\ref{thm:BBP2}
also hold with $\lambda_{\max}$ replaced by $\frac{L(N,M)}{M}$. In
particular, Theorem~\ref{thm:BBP1} shows that as long as the site
passage time on the distinguished columns have mean less than
$1+\gamma^{-1}$, the last passage time has the same limit behavior
as the case when all the sites are identically distributed.

\subsection{Queues in tandem}

Suppose that there are $N$ servers and $M$ customers. Initially
all the customers are at the first server in a queue. Once a
customer is served at a server, then (s)he moves to the queue of
the next server and waits for his/her turn. The service time for
the $j$th customer at the $i$th server is assumed to be a random
variable $X(i,j)$ and let $D(N,M)$ be the departure time of all
the customers from all the queues. It is well-known that $D(N,M)$
has the same distribution as $L(N,M)$ of the last passage
percolation model (see e.g. \cite{GlynnWhitt}).

In the queueing theory context, Theorem~\ref{thm:BBP1} determines
the effect of a few slow servers to the total departure time.
Suppose that $X(i,j)$ is an exponential random variable of mean
$1$ for $i=r+1, \dots, N$ (independent of $i$) and of mean
$\ell_i$ for $i=1, \dots, r$. In other words, the service times at
the first $r$ servers are distributed differently from those at
the rest of the servers. When all of $\ell_i$ are not so large,
the departure time has the same limiting law as when all the
service times are identically distributed, but the whole process
slows down when some of the servers are sufficiently slow.
Theorem~\ref{thm:BBP1} shows that the critical value is
$\ell_i=1+\gamma^{-1}$. Note that due to a symmetry between
servers and customers, the theorem also applies to slow customers.

\subsection{Totally asymmetric simple exclusion process}

The last passage percolation can also be interpreted as an
interacting particle systems (see e.g. \cite{Seppal1997,
kurtj:shape}). We will consider the totally asymmetric simple
exclusion process. Let $x_j(t)\in \mathbb{Z}$, $x_j(t)$,
$j=1,2,\dots$, $t\in [0,\infty)$, denote the location of the $j$th
particle at time $t$. A particle can jump only to its right
neighboring site after random time if the site is not occupied.
Let $X(i,j)$ be independent random variables which represent the
$i$th jumping time of the $j$th particle. We take the initial
condition as $x_j(0)=1-j$, $j=1,2,\dots$. Then $X(i,j)$ is the
time it takes for the $j$th particle $x_j$ to jump from the site
$i-j$ to $i-j+1$.

Let $T(i,j)$ be the time it takes for the $j$th particle to arrive
at the location $i-j+1$. Equivalently, $T(i,j)$ is the time it
takes for the $j$th particle to perform the first $i$ jumps. Note
that in order for the $j$th particle to jump from the site $i-j$
to $i-j+1$, the $(j-1)$th particle should be to the right of the
site $i-j+1$. Hence we find that
\begin{equation}
  T(i,j)= \max\{ T(i-1,j), T(i,j-1) \} + X(i,j), \qquad i,j\ge 1,
\end{equation}
where $T(0,j)=T(i,0)=0$, $i,j\ge 1$, by definition. A simple
geometric consideration shows that last passage time $L(i,j)$
satisfies exactly the same recurrence relation. Therefore $T(i,j)$
is same as the last passage time in the sense of distribution.

Let $\#(m,t)$ denote the number of particles to the right of the
site $m$ at time $t$. The flux $F(m,t)$, the number of particles
that have jumped cross the interval $(m,m+1)$ up to time $t$, is
then $F(m,t)=\#(m,t)$ for $m>0$, and $F(m,t)=\#(m,t)+m$ for $m\le
0$. The event that $\#(m,t)\ge M$ is same as the event that the
$M$th particle is to the right of the site $m$ at time $t$. This
is again equal to the event that $T(m+M, M)\le t$, and hence we
find that $\mathbb{P}(\#(m,t)\ge M)= \mathbb{P}( T(m+M, M)\le t) =
\mathbb{P}(L(m+M, M)\le t)$. Therefore, Theorem~\ref{thm:BBP1} and
Theorem~\ref{thm:BBP2} again apply to $\#([ut],t)$, and hence to
$F([ut], t)$. We state the results for $\#([ut], t)$ here.

\subsubsection*{Traffic of slow start from stop}

Suppose that $X(i,j)$ is an independent exponential random
variable of mean $\ell_i$ for $i=1, \dots, r$ and of mean $1$ for
$i>r$ (independent of $j$). In other words, each particle jumps at
rate $\frac1{\ell_i}$ for its first $r$ jumps and then jumps at
rate $1$ afterwards. When $\ell_1\ge \dots\ge \ell_r$, one can
view it as a toy model for the following traffic situation:
(infinite) cars in one lane, which were fully stopped at the red
signal, speed up at the green signal until they finally reach the
steady speed (after $r$ `jumps'). Set $\ell=\max\{\ell_1, \dots,
\ell_r\}$ and let $k\ge 1$ be the number of $\ell_i$'s equal to
$\ell$. By re-interpreting Theorem~\ref{thm:BBP1}, a tedious but
straightforward calculation shows the following results for $-1<
u\le 0$ :
\begin{equation}\label{eq:slowtraffic1}
  \lim_{t\to\infty} \mathbb{P}\biggl( \#([ut], t) \ge \frac14(1-u)^2 t + x \biggl(
  \frac{1-u^2}{4} \biggr)^{2/3} t^{1/3} \biggr)
  = \begin{cases} F_0(-x), \qquad & u\in \bigl( 1-\frac2{\ell}, 0]\cap (-1, 0] \\
  F_k(-x), \qquad &u=1-\frac2{\ell} \in (-1, 0],
  \end{cases}
\end{equation}
and for $u\in  \bigl(\frac1{-\ell}, 1-\frac2{\ell} \bigr) \cap
(-1,0]$,
\begin{equation}\label{eq:ubig}
  \lim_{t\to\infty} \mathbb{P}\biggl( \#([ut], t) \ge \frac{\ell-1-\ell u}{\ell^2} t + x
  \frac{(\ell-1)^{3/2}(\ell-1-\ell u)}{\ell^{9/2}(\ell-2-\ell u)^{1/2}} t^{1/2} \biggr)
  = G_k(-x).
\end{equation}
This shows that fast jumps do not affect the flux but slow jumps
may change the flux. When $r=0$ (all cars jumping at the same
rate), ~\eqref{eq:slowtraffic1} was first obtained in
\cite{kurtj:shape} for $0\le u<1$.

\subsubsection*{Traffic with a few slow cars}

The exclusion process of the particles yields a dual process of
the holes. As the particles jump to the right, the holes, the
unoccupied sites, jump to the left. The leftmost hole jumps at
rate $\frac1{\ell_1}$ since each particle jump at that rate at its
first jump. Likewise, the second leftmost hole jumps at rate
$\frac1{\ell_2}$ and so on. Hence this model can be thought of
traffic model where there are a few cars of distinguished jump
rates. Initially the holes are at the sites $\{1,2,3,\dots\}$. As
the number of holes $\mathcal{H}(m,t)$ on the left of the site
$m+1$ at time $t$ satisfies $\mathcal{H}(m,t)=\#(m,t)+m$,
~\eqref{eq:slowtraffic1} and~\eqref{eq:ubig} imply the following
results for $-1< u\le 0$ :
\begin{equation}
  \lim_{t\to\infty} \mathbb{P}\biggl( \mathcal{H}([ut], t) \ge \frac14(1+u)^2 t + x \biggl(
  \frac{1-u^2}{4} \biggr)^{2/3} t^{1/3} \biggr)
  = \begin{cases} F_0(-x), \qquad & u\in \bigl( 1-\frac2{\ell}, 0]\cap (-1, 0] \\
  F_k(-x), \qquad &u=1-\frac2{\ell} \in (-1, 0],
  \end{cases}
\end{equation}
and for $u\in  \bigl(\frac1{-\ell}, 1-\frac2{\ell} \bigr) \cap
(-1,0]$,
\begin{equation}
  \lim_{t\to\infty} \mathbb{P}\biggl( \mathcal{H}([ut], t) \ge \frac{\ell-1-\ell u+\ell^2 u}{\ell^2} t + x
  \frac{(\ell-1)^{3/2}(\ell-1-\ell u)}{\ell^{9/2}(\ell-2-\ell u)^{1/2}} t^{1/2} \biggr)
  = G_k(-x).
\end{equation}

\bigskip

The full case of $-1<u<1$ and also correlation functions of
various locations for both of the above traffic models will be
discussed in a forthcoming paper.

\section{Proof of Theorem \ref{thm:det}}\label{sec:proofthmdet}

We prove Theorem~\ref{thm:det} in this section.

Since both sides of \eqref{eq:detformula} are analytic in each
$w_j$, the case when some of $w_j$'s coincide follows from
analytic continuation of the case when all $w_j$'s are distinct.
Hence we assume in this section that all $w_j$'s are distinct. We
need to prove that
\begin{equation}\label{eq:detiden}
\begin{split}
  &\det \biggl( \delta_{mn} -
< \frac1{1-\mathbf{A}_x} s^{(m)}(w_1,\dots, w_m),
t^{(n)}(w_1,\dots, w_{n-1}) >_{(x,\infty)}
  \biggr)_{1\le m,n \le k} \\
 &  =  \frac{1}{\displaystyle \prod_{1\le m<n\le k} (w_n-w_m)}
  \det \begin{pmatrix}
  (w_m+ D_x)^{n-1}f(x,w_m)
  \end{pmatrix}_{1\le m,n\le k}.
\end{split}
\end{equation}

\begin{notrem}
In the below, we sometimes have a product of empty indices. For
instance when $n=1$, the product $\displaystyle \prod_{a=1}^{n-1}
(w_a-w_i)$ in \eqref{eq:EinS} has no indices. In such cases, we
interpret the product as $1$.
\end{notrem}

\bigskip

Let $\mathbb{A}_x : L^2((0,\infty)) \to L^2((0,\infty))$ be the operator
with kernel
\begin{equation}
  \mathbb{A}_x(u,v)=\mathbf{A}(u+x, v+x).
\end{equation}
Set
\begin{equation}\label{eq:defST1}
  S_m(u)= s^{(m)}(u+x) = s^{(m)}(u+x; w_1, \dots, w_m)
\end{equation}
and set
\begin{equation}\label{eq:defST2}
  T_m(v)=t^{(m)}(v+x) = t^{(m)}(v+x; w_1, \dots, w_{m-1}).
\end{equation}
Then
\begin{equation}
   \biggl( \frac1{1-\mathbf{A}_x} s^{(m)} \biggr) (u+x)
=  \biggl(\frac1{1-\mathbb{A}_x} S_m \biggr)(u)
\end{equation}
and the matrix on the left-hand-side of \eqref{eq:detiden} is
\begin{equation}
  \bigl( \delta_{ij} - < \frac1{1-\mathbb{A}_x} S_i, T_j>\bigr)_{1\le i,j\le k}
\end{equation}
where $<,>=<,>_{(0,\infty)}$ is the real inner product in $L^2((0,\infty))$.

Since
\begin{equation}
  \prod_{j=1}^m \frac{1}{w_j+ia} = \sum_{j=1}^m \frac1{w_j+ia}
  \displaystyle \prod_{\substack{\ell=1 \\ \ell\neq
  j}}^m \frac1{w_\ell-w_j},
\end{equation}
we find
\begin{equation}\label{eq:SinE}
   \biggl( \frac1{1-\mathbb{A}_x} S_m \biggr) (u)
   = \sum_{j=1}^m
  \biggl[ \displaystyle \prod_{\substack{\ell=1 \\ \ell\neq
  j}}^m \frac1{w_\ell-w_j} \biggr]
  E_{w_j}(u).
\end{equation}
where
\begin{equation}
  E_w(u)=E_{w}(u;x) :=
  \biggl( \frac1{1-\mathbb{A}_x} \widetilde{C}_{w} \biggr) (u),
\qquad \widetilde{C}_w(u) := C_w (u+x)
\end{equation}
(recall~\eqref{eq:Cdef} for the definition of $C_w$). For later
use, we note that $f(x,w)$ defined in \eqref{eq:fdef} satisfies
that
\begin{equation}\label{eq:recallf}
  f(x,w) = 1- <E_w, T_1>.
\end{equation}
Now we invert the relation \eqref{eq:SinE}. For $1\le m\le k$,
\eqref{eq:SinE} is a system of $k$ linear equations for $E_{w_j}$,
$1\le j\le k$.

\begin{lem}
The equation \eqref{eq:SinE} for $E_{w_j}$ has the solution given
by
\begin{equation}\label{eq:EinS}
  E_{w_j}(u)  = \sum_{n=1}^j \biggl[ \prod_{a=1}^{n-1} (w_a-w_j)
  \biggr]
   \biggl( \frac1{1-\mathbb{A}_x} S_n \biggr) (u).
\end{equation}
\end{lem}

\begin{proof}
Consider the function
\begin{equation}
  F(z) := -\prod_{\ell =n}^m \frac1{w_\ell-z}.
\end{equation}
Integrating over a circle of radius $R$, and then taking
$R\to\infty$, we find that the sum of residues of $F$ is equal to
$0$ when $m-n\ge 1$ and is equal to $1$ when $m=n$. On the other
hand, by directly computation, the residue of $F$ at $z=w_j$ is
\begin{equation}
  \prod_{\substack{ \ell=n \\ \ell \neq j}}^m \frac1{w_\ell -w_j}.
\end{equation}
Hence we obtain the identity
\begin{equation}\label{eq:deltaiden}
  \sum_{j=n}^m \prod_{\substack{ \ell=n \\ \ell \neq j}}^m \frac1{w_\ell
  -w_j} = \delta_{mn}, \qquad m\ge n.
\end{equation}

Now as all $w_i$'s are distinct, the determinant of the matrix for
the linear equation \eqref{eq:SinE} is $\prod_{1\le \ell<m\le
k}(w_{\ell}-w_m)^{-1}$, which is non-zero. Hence there is a unique
solution $E_{w_j}$ for \eqref{eq:SinE}. We should check that
\eqref{eq:EinS} solves \eqref{eq:SinE}. But this follows by
inserting \eqref{eq:EinS} into the right-hand-side of
\eqref{eq:SinE}, changing the order of summations, and then using
\eqref{eq:deltaiden}.
\end{proof}

From \eqref{eq:EinS}, we obtain for each $1\le i,j\le k$,
\begin{equation}\label{eq:rowoppp}
  \sum_{n=1}^i \prod_{a=1}^{n-1} (w_a-w_i) \cdot \bigl(\delta_{nj}
- < \frac1{1-\mathbb{A}_x} S_n, T_j> \bigr) = F_{ij} - <E_{w_i}, T_j>
\end{equation}
where
\begin{equation}\label{eq:Finitial}
  F_{ij} :=
  \prod_{a=1}^{j-1} (w_a-w_i).
\end{equation}
Note that $F_{ij}=0$ when $i<j$. Now we perform row operations of
the matrix $\bigl(\delta_{ij}-<\frac1{1-\mathbb{A}_x} S_i,
T_j>\bigr)_{1\le i,j\le k}$ using~\eqref{eq:rowoppp} that replaces
the $i$th row by a linear combination of the first $i$ rows to
find that
\begin{equation}\label{eq:detid2}
\begin{split}
  \det\bigl(\delta_{ij} - <\frac1{1-\mathbb{A}_x} S_i, T_j>
\bigr)_{k\times k}
  &= \prod_{i=1}^k\prod_{a=1}^{i-1} \frac1{w_a-w_i}
  \cdot \det \bigl( F_{ij} - <E_{w_i}, T_j > \bigr)_{k\times k} \\
  &= \frac1{\displaystyle \prod_{1\le i<j\le k} (w_i-w_j)} \cdot
  \det \bigl( F_{ij} - <E_{w_i}, T_j > \bigr)_{k\times k} .
\end{split}
\end{equation}
Note that the when $j=1$, $F_{ij}=1$, (see the Notational Remark
above) and hence the first column of the matrix $\bigl( F_{ij} -
<E_{w_i}, T_j > \bigr)_{k\times k}$ consists of the functions (see
\eqref{eq:recallf})
\begin{equation}
   1-<E_{w_i}, T_1>= f(x,w_i).
\end{equation}
For example, when $k=3$, the determinant on the right-hand-side of
\eqref{eq:detid2} is
\begin{equation}
  \det \begin{pmatrix} f(x,w_1) & -<E_{w_1}, T_2> & -<E_{w_1}, T_3> \\
  f(x,w_2) & (w_1-w_2) -<E_{w_2}, T_2> & -<E_{w_2}, T_3 > \\
  f(x,w_3) & (w_1-w_3) -<E_{w_3}, T_2> & (w_1-w_3)(w_2-w_3) -<E_{w_3}, T_3>
  \end{pmatrix}.
\end{equation}

From the definition \eqref{eq:defST2} of $T_j$ and the definition
\eqref{eq:tdef} of $t_j$, we have
\begin{equation}\label{eq:trel}
  T_j = w_{j-1}T_{j-1} - D_xT_{j-1}, \qquad j\ge 2.
\end{equation}
Set $M^{(0)}$ be the matrix
\begin{equation}\label{eq:M0def}
  M^{(0)} := \bigl( F_{ij}-<E_{w_i},T_j> \bigr)_{1\le i,j\le k}.
\end{equation}
Let $M^{(1)}$ be the matrix defined by
\begin{equation}
  M^{(1)} := M^{(0)} \begin{pmatrix}
  1 & -w_1 & 0 \\
  0 & 1 & -w_2 & 0 \\
    & 0 & 1 & -w_3  & 0\\
  & & & \ddots & \ddots & \\
  & & & &  \ddots & \ddots  &\\
  & & & & 0 & 1 & -w_{k-1} \\
  & & & &   & 0 & 1
  \end{pmatrix},
\end{equation}
whose determinant is same as the determinant of $M^{(0)}$. Using
the relation \eqref{eq:trel}, the entries of $M^{(1)}= \bigl(
M^{(1)}_{ij} \bigr)_{1\le i,j\le k}$ are given by
\begin{equation}
  M^{(1)}_{ij} =
  \begin{cases}
  f(x,w_i), \qquad & j=1 \\
  - \bigl( F_{ij}^{(1)} - <E_{w_i}, D_x T_{j-1} >), \qquad & 2\le j\le
  k,
  \end{cases}
\end{equation}
where
\begin{equation}
  F_{ij}^{(1)} := w_{j-1}F_{i,j-1}- F_{ij}.
\end{equation}

Now define a new matrix $M^{(2)}=\bigl( M^{(2)}_{ij} \bigr)_{1\le
i,j\le k}$ by
\begin{equation}
  M^{(2)} := M^{(1)} \begin{pmatrix}
  1 & 0 & 0 \\
  0 & 1 & -w_1 & 0 \\
    & 0 & 1 & -w_2  & 0\\
  & & & \ddots & \ddots & \\
  & & & &  \ddots & \ddots  &\\
  & & & & 0 & 1 & -w_{k-2} \\
  & & & &   & 0 & 1
  \end{pmatrix}.
\end{equation}
Using the relation
\begin{equation}
  D_xT_{j-1} = w_{j-2}D_xT_{j-2} - D_x^2T_{j-2}
\end{equation}
that follows from \eqref{eq:trel} for $3\le j\le k$, we find that
\begin{equation}
  M^{(2)}_{ij} =
  \begin{cases}
  f(x,w_i), \qquad &j=1 \\
  - \bigl( F_{i2}^{(1)} - <E_{w_i}, D_x T_{1} >), \qquad &j=2 \\
  F_{ij}^{(2)} - <E_{w_i}, D_x^2 T_{j-2} >, \qquad & 3\le j\le
  k,
  \end{cases}
\end{equation}
where
\begin{equation}
  F_{ij}^{(2)} :=  w_{j-2}F_{i,j-1}^{(1)} -F_{ij}^{(1)}.
\end{equation}

Continuing in a similar way, we eventually define
\begin{equation}
\begin{split}
  M^{(k-1)} := M^{(0)} & \begin{pmatrix}
  1 & -w_1 & 0 \\
  0 & 1 & -w_2 & 0 \\
    & 0 & 1 & -w_3  & 0\\
  & & & \ddots & \ddots & \\
  & & & &  \ddots & \ddots  &\\
  & & & & 0 & 1 & -w_{k-1} \\
  & & & &   & 0 & 1
  \end{pmatrix}
  \times
  \begin{pmatrix}
  1 & 0 & 0 \\
  0 & 1 & -w_1 & 0 \\
    & 0 & 1 & -w_2  & 0\\
  & & & \ddots & \ddots & \\
  & & & &  \ddots & \ddots  &\\
  & & & & 0 & 1 & -w_{k-2} \\
  & & & &   & 0 & 1
  \end{pmatrix}
  \\
  &
  \times
  \cdots
  \times
  \begin{pmatrix}
  1 & 0 & 0 \\
  0 & 1 & 0 & 0 \\
    & 0 & 1 & 0  & 0\\
  & & & \ddots & \ddots & \\
  & & & &  \ddots & \ddots  &\\
  & & & & 0 & 1 & -w_{1} \\
  & & & &   & 0 & 1
  \end{pmatrix}.
\end{split}
\end{equation}
Using the fact that for all $\ell \ge 1$, $j\ge 2$,
\begin{equation}
  D_x^{\ell} T_{j} = w_{j-1}D_x^{\ell}T_{j-1} - D_x^{\ell +1}
  T_{j-1},
\end{equation}
we find that
\begin{equation}\label{eq:Mkdef}
  M^{(k-1)}_{ij} =
  (-1)^{j-1} \bigl( F_{ij}^{(j-1)} - <E_{w_i}, D_x^{j-1} T_{1} >
  \bigr), \qquad 1\le i, j\le k,
\end{equation}
where for each $1\le i\le k$, $F_{ij}^{(\ell)}$ is inductively
defined by the relation
\begin{equation}\label{eq:rec}
  F_{ij}^{(\ell)} :=  w_{j-\ell}F_{i,j-1}^{(\ell-1)} -F_{ij}^{(\ell-1)} ,
  \qquad 1\le \ell < j\le k
\end{equation}
and (see \eqref{eq:Finitial})
\begin{equation}\label{eq:Fzero}
    F_{ij}^{(0)} := F_{ij} = \prod_{a=1}^{j-1} (w_a-w_i), \qquad 1\le j\le k.
\end{equation}
Recall that $F_{i1}=1$, and hence $M^{(k-1)}_{i1}=1-<E_{w_i}, T_1>
= f(x,w_i)$.

\begin{lem}
The solution $F_{ij}^{(\ell)}$ to the recurrence relation
\eqref{eq:rec} and \eqref{eq:Fzero} is
\begin{equation}\label{eq:recsol}
  F_{ij}^{(\ell)} = w_i^{\ell} \prod_{a=1}^{j-\ell-1}(w_a-w_i),
  \qquad 1\le \ell <j\le k.
\end{equation}
\end{lem}

\begin{proof}
This follows easily from an induction in $\ell$. Here, as
mentioned in the Notational Remark above, when $j=\ell+1$, we
understand that the product $\prod_{a=1}^{0} (w_a-w_i)=1$.
\end{proof}

Therefore $F^{(j-1)}_{ij}=w_i^{j-1}$, and as $\det(M^{(0)})=
\det(M^{(k-1)})$, we find from \eqref{eq:detid2}, \eqref{eq:M0def}
and \eqref{eq:Mkdef} that
\begin{equation}\label{eq:middle}
\begin{split}
  \det\bigl(\delta_{ij} - <\frac1{1-\mathbb{A}_x} S_i, T_j>\bigr)_{1\le i,j\le k}
  = \frac{(-1)^{[k/2]}}{\displaystyle \prod_{1\le i<j\le k} (w_i-w_j)}
  \det \bigl( M \bigr),
\end{split}
\end{equation}
where $[k/2]$ denotes the largest integer smaller than or equal to
$k/2$, and the $k\times k$ matrix $M= \bigl( M_{ij}\bigr)_{1\le
i,j\le k}$ is given by
\begin{equation}
  M_{ij} =
  w_i^{j-1} - <E_{w_i}, D_x^{j-1} T_{1} >.
\end{equation}
As $\prod_{1\le m<n\le k} (-1)= (-1)^{[k/2]}$, this is equal to
\begin{equation}\label{eq:middle2}
\begin{split}
  \det\bigl(\delta_{ij} - <\frac1{1-\mathbb{A}_x} S_i, T_j>\bigr)_{1\le i,j\le k}
  = \frac{1}{\displaystyle \prod_{1\le i<j\le k} (w_j-w_i)}
  \det \bigl( M \bigr),
\end{split}
\end{equation}

Now we will show that $\det(M)$ is equal to the determinant on the
right-hand-side of \eqref{eq:detiden}. For this purpose, we use
the following result.

\begin{lem}
For $\ell\ge 0$, there are smooth functions $F_{\ell, a}(x)$,
$a=0,1,\dots, \ell-1$, such that
\begin{equation}\label{eq:wDgoal}
  (w + D_x)^{\ell} f(x,w) = w^{\ell} - <E_w, D_x^\ell T_1> -
  \sum_{a=0}^{\ell-1} F_{\ell, a}(x) (w+D_x)^a f(x,w) .
\end{equation}
\end{lem}

\begin{proof}
From the definition of $E_w$,
\begin{equation}
  D_x E_w = D_x \biggl( \frac1{1-\mathbb{A}_x}\widetilde{C}_w \biggr)
  = \frac{1}{1-\mathbb{A}_x} (D_x\mathbb{A}_x) \frac1{1-\mathbb{A}_x} \widetilde{C}_w +
  \frac1{1-\mathbb{A}_x}D_x\widetilde{C}_w.
\end{equation}
It is direct to check that $(D_x\mathbb{A}_x)(u,v) =
-\Ai(x+u)\Ai(x+v)$. Hence $D_x\mathbb{A}_x=-T_1\otimes T_1$. On
the other hand, from the definition \eqref{eq:Cdef} of $C_w$,
$D_x\widetilde{C}_w= T_1-w\widetilde{C}_w$. Hence we find
\begin{equation}
\begin{split}
  D_xE_w
  &= -\frac{1}{1-\mathbb{A}_x} T_1\otimes T_1 \frac1{1-\mathbb{A}_x} \widetilde{C}_w +
  \frac1{1-\mathbb{A}_x}(T_1-w\widetilde{C}_w) \\
  &= - < T_1, \frac1{1-\mathbb{A}_x} \widetilde{C}_w> \frac{1}{1-\mathbb{A}_x} T_1
  + \frac1{1-\mathbb{A}_x}T_1 -w \frac1{1-\mathbb{A}_x}\widetilde{C}_w,
\end{split}
\end{equation}
which implies that
\begin{equation}\label{eq:wDE}
  (w+D_w)E_w = f(x,w)\frac{1}{1-\mathbb{A}_x} T_1
\end{equation}

Now we use an induction in $\ell$ to prove \eqref{eq:wDgoal}. When
$\ell=0$, by definition \eqref{eq:fdef} of $f$, \eqref{eq:wDgoal}
holds. Now suppose that \eqref{eq:wDgoal} holds true for some
$\ell\ge 0$. Then using the general identities $(w+D_x)<h,g>=
<(w+D_x)h, g> + <h,D_xg>$ and $(w+D_x) \bigl(hg) = (D_xh)g+ h
\bigl((w+D_x)g \bigr)$,
\begin{equation}
\begin{split}
  & (w + D_x)^{\ell+1} f(x,w) \\
  &\qquad = (w + D_x) \bigl[ (w + D_x)^{\ell} f(x,w) \bigr] \\
  &\qquad = (w + D_x) \bigl[ w^{\ell} - < E_w, D_x^{\ell} T_1> -
  \sum_{a=0}^{\ell-1} F_{\ell, a}(x) (w+D_x)^a f(x,w) \bigr]\\
  &\qquad  = w^{\ell+1} - <(w+D_x)E_w, D_x^\ell T_1> -  <E_w, D_x^{\ell+1}
  T_1> \\
  &\qquad\quad  - \sum_{a=0}^{\ell-1} \biggl\{
  \bigl( D_xF_{\ell, a}(x)\bigr)  (w+D_x)^{a} f(x,w)
  +  F_{\ell, a}(x) (w+D_x)^{a+1} f(x,w)  \biggr\} \\
  &\qquad  = w^{\ell+1} - f(x,w) < \frac1{1-\mathbb{A}_x}T_1, D_x^\ell T_1>
  -  <E_w, D_x^{\ell+1} T_1>  \\
  &\qquad\quad  - \sum_{a=0}^{\ell-1} \biggl\{
  \bigl( D_xF_{\ell, a}(x)\bigr)  (w+D_x)^{a} f(x,w)
  +  F_{\ell, a}(x) (w+D_x)^{a+1} f(x,w)  \biggr\} ,
\end{split}
\end{equation}
where \eqref{eq:wDE} is applied in the last step. Therefore we
find that \eqref{eq:wDgoal} holds true for $\ell+1$ with the
functions
\begin{equation}
  F_{\ell+1, a} =
  \begin{cases}
   D_xF_{\ell, 0}(x) + <\frac1{1-\mathbb{A}_x}T_1, D_x^\ell T_1> , \qquad & a=0, \\
   D_xF_{\ell, a}(x) + F_{\ell, a-1}, \qquad & 1\le a\le \ell-1,
   \\
   F_{\ell, \ell-1}, \qquad & a=\ell,
  \end{cases}
\end{equation}
where $F_{0,0}:=0$.
\end{proof}

By applying \eqref{eq:wDgoal} repeatedly to $(w+D_x)^af(x;w)$
inside the summation on the right-hand-side of \eqref{eq:wDgoal},
for $\ell\ge 0$, there are smooth functions $G_{\ell, a}$,
$a=0,1,\dots, \ell-1$ such that
\begin{equation}\label{eq:wE}
  w^{\ell} - <E_w, D_x^\ell T_1>
  =
  (w + D_x)^{\ell} f(x;w) +
  \sum_{a=0}^{\ell-1} G_{\ell, a}(x) \bigl( w^{a} - <E_w, D_x^a T_1> \bigr).
\end{equation}
Therefore for any $1\le i\le k$,
\begin{equation}
  M_{ij} =
  w_i^{j-1} - <E_{w_i}, D_x^{j-1} T_{1} >
  = (w_i+D_x)^{j-1} f(x,w_i) + \sum_{a=1}^{j-1} G_{j-1, a-1}(x)
  M_{ia} .
\end{equation}
In other words, the $j$th column vector in the matrix $M$ is
equal to a linear combination of the first, second, ..., $j-1$th
column vectors plus the vector $\bigl( (w_1+D_x)^{j-1} f(x,w_1),
\cdots, (w_k+D_x)^{j-1} f(x,w_k) \bigr)^T$. Hence by applying
proper column operations, we find
\begin{equation}
  \det(M) = \det\bigl( (w_i+D_x)^{j-1} f(x,w_i) \bigr)_{1\le i,j\le
  k}.
\end{equation}
This, together with \eqref{eq:middle2}, implies that the
left-hand-side of \eqref{eq:detiden} is equal to
\begin{equation}
\begin{split}
  \frac{1}{\displaystyle \prod_{1\le m<n\le k} (w_n-w_m)}
  \det \begin{pmatrix}
  (w_m+ D_x)^{n-1}f(x,w_m)
  \end{pmatrix}_{1\le m,n\le k}.
\end{split}
\end{equation}
This is the right-hand-side of \eqref{eq:detiden} and
Theorem~\ref{thm:det} is proved.


\section{Proof Proposition~\ref{prop:pidet}}
\label{sec:ProofPropdidet}

Let $T_n(\varphi)=(\varphi_{i-j})_{0\le i,j\le n-1}$ be the
Toeplitz matrix of the symbol $\varphi$ where $\varphi_k$ denotes
the Fourier coefficients of $\varphi$. Let $D_n(\varphi)=\det
T_n(\varphi)$ be the Toeplitz determinant. We will use the
following identity \cite{GeronimoC69, BoOk} between a Toeplitz
determinant and the Fredholm determinant of an operator on a
discrete set: for all $n\ge 1$,
\begin{equation}\label{eq:GCBO}
  \frac{D_n(\varphi)}{G(\varphi)^nE(\varphi)} = \det (1- P_n AB P_n)
\end{equation}
with
\begin{equation}\label{eq:Edef}
  G(\varphi)=e^{ (\log \varphi)_0 },
  \qquad E(\varphi)
  = \exp\biggl\{\sum_{k=1}^\infty k(\log \varphi)_k (\log
  \varphi)_{-k}\biggr\}
\end{equation}
where the operators $A, B$ are defined by the kernels
\eqref{eq:defAB}. This identity holds, for example, for
complex-valued analytic functions $\varphi$ with zero winding
number which has a Weiner-Hopf factorization. See
\cite{BottSimple} for the minimal condition on $\varphi$ for which
the identity holds.

It is well-known that $\pi_n^*$ has the multi-integral expression
(see e.g. \cite{Szego})
\begin{equation}\label{eq:mltint}
  \pi_n^*(z) = \frac1{D_n(\phi)} \int_{|z_1|=1} \dots \int_{|z_n|=1}
  \prod_{1\le j<k\le n} (1-z_kz_j^{-1})  \prod_{j=1}^n (1-z z_j^{-1})
  \prod_{j=1}^n \phi(z_j) \frac{dz_j}{2\pi i z_j} .
\end{equation}
By using the multi-integral formula of a Toeplitz determinant,
\eqref{eq:mltint} can be written as
\begin{equation}\label{eq:piDD}
  \pi_n^*(z) = \frac{D_n(\phi_z)}{D_n(\phi)}
\end{equation}
where the new symbol $\phi_z$ is
\begin{equation}
  \phi_z(w) := \biggl(1-\frac{z}{w}\biggr) \phi(w).
\end{equation}
See \cite{BDS} for a use of the identity~\eqref{eq:piDD} in random
matrix theory. Using~\eqref{eq:GCBO} for $D_n(\phi_z)$ and
$D_n(\phi)$, we find that
\begin{equation}
\begin{split}
  \pi_n^*(z) &= \frac{G(\phi_z)^nE(\phi_z)}{G(\phi)^nE(\phi)} \cdot
  \frac{\det(1-P_n A_z B_z P_n)}{\det(1-P_n A B P_n)}\\
  &= \frac{G(\phi_z)^nE(\phi_z)}{G(\phi)^nE(\phi)} \cdot \det\bigg(
  1- \frac1{1-P_n A B P_n} P_n (A_zB_z-AB)P_n \biggr)
\end{split}
\end{equation}
where the operators $A_z, B_z : \ell^2(\mathbb{N}_0) \to
\ell^2(\mathbb{N}_0)$ are defined by the kernels
\begin{equation}
  A_z(j,m) = (\psi_z^{-1})_{j+m+1}, \qquad
  B_z(m,k) = (\psi_z)_{-m-k-1}.
\end{equation}
where $\psi_z=(\phi_z)_+/(\phi_z)_-$.

When $|z|<1$, from \eqref{eq:Edef}, it is easy to check that
\begin{equation}
  \frac{G(\phi_z)^nE(\phi_z)}{G(\phi)^nE(\phi)}
  = e^{-\sum_{k=1}^\infty (\log\phi)_{k} z^k}. 
\end{equation}
Now we consider $A_zB_z-AB$. As $\phi_z$ has the Wiener-Hopf
factorization $\phi_z=(\phi_z)_+(\phi_z)_-$ where $(\phi_z)_+(w) =
\phi_+(w)$ and $(\phi_z)_-(w)=\bigl(1-\frac{z}{w}\bigr)\phi_-(w)$,
we find $\psi_z(w)= \frac1{(1-\frac{z}{w})}\psi(w)$. Therefore,
\begin{equation}
\begin{split}
  (A_zB_z)(j,k)-(AB)(j,k)
  &= \sum_{m=0}^\infty (\psi_z^{-1})_{j+m+1} (\psi_z)_{-m-k-1}
  - (\psi^{-1})_{j+m+1} \psi_{-m-k-1} \\
  & = \sum_{m=0}^\infty
  \int_{|a|=1} \int_{|b|=1} a^{-j-m-1}b^{m+k+1}
  \biggl( \frac{1-\frac{z}{a}}{1-\frac{z}{b}}-1 \biggr)
  \frac{\psi(b)}{\psi(a)} \frac{da}{2\pi ia} \frac{db}{2\pi ib}.
\end{split}
\end{equation}
Note that $\phi_+(z)=\phi(z)/\phi_-(z)$ is analytic in a region
outside the unit circle as $\phi(z)$ is assumed to be analytic in
a neighborhood of the unit circle and $\phi_-(z)$ is analytic
outside the unit circle. Hence $\phi_+(z)$ is analytic in a
neighborhood of the unit circle. Similarly, $\phi_-(z)$ is
analytic in a neighborhood of the unit circle and so is $\psi(z)$.
Therefore the contours $|a|=1$ and $|b|=1$ can be deformed so that
$|a|>|b|$. Hence,
\begin{equation}
\begin{split}
  (A_zB_z)(j,k)-(AB)(j,k)
  &= \int\int_{|a|>|b|} a^{-j-1}b^{k+1} \biggl[ \sum_{m=0}^\infty
  \biggl( \frac{b}{a} \biggr)^m \biggr]
  \frac{z(a-b)}{a(b-z)}
  \frac{\psi(b)}{\psi(a)} \frac{da db}{(2\pi i)^2ab} \\
  &= \int\int_{|a|>|b|} a^{-j-1}b^{k+1}
  \frac{z}{b-z}
  \frac{\psi(b)}{\psi(a)} \frac{da db}{(2\pi i)^2ab}
  = Q(j)R(k)
\end{split}
\end{equation}
where $Q$ and $R$ are defined by \eqref{eq:defQR}. This implies
that $A_zB_z$ is a rank 1 perturbation of $AB$
and we find
\begin{equation}
\begin{split}
  \pi_n^*(z)
  &= e^{-\sum_{k=1}^\infty (\log\phi)_{k} z^k} \cdot \det\biggl(
  1- \frac1{1-P_n A B P_n} P_n Q \otimes R P_n \biggr) \\
  &= e^{-\sum_{k=1}^\infty (\log\phi)_{k} z^k} \biggl\{ 1- <\frac1{1-P_n A B P_n} P_n Q,
  P_nR> \biggr\}
\end{split}
\end{equation}
which completes the proof of \eqref{eq:pidet}.

Proof for \eqref{eq:pidet222} is similar by noting that
\begin{equation}
    \pi_n(z)= \frac{D_n(\phi^z)}{D_n(\phi)}, \qquad
    \phi^z(w):=(z-w)\phi(w).
\end{equation}

\section{Proof of
Theorem~\ref{thm:fPain}}\label{sec:proofTheoremfPain}

We apply Proposition \ref{prop:pidet} to the function
\begin{equation}\label{eq:symbolsp}
  \phi(z):= e^{t(z+z^{-1})}
\end{equation}
for positive number $t$. Then \eqref{eq:pidet} becomes
\begin{equation}\label{eq:piid}
  e^{tz} \pi_n^*(z)= 1- <\frac1{1-P_n A B P_n}
  P_nQ, P_nR>_{\ell^2(\mathbb{N}_0)}
\end{equation}
with
\begin{equation}\label{eq:5.3}
  \psi(z) = e^{t(z-\frac1{z})}.
\end{equation}
It is easy to check that the inner product on the right-hand-side
is unchanged when the functions $A(j,m)$, $B(m,k)$, $Q(j)$ and
$R(k)$ are replaced by
\begin{equation}\label{eq:signchange}
  (-1)^{j+m}A(j,m), \quad (-1)^{m+k}B(m,k), \quad (-1)^jQ(j),
  \quad (-1)^k R(k),
\end{equation}
respectively. We will denote these new functions by the same
notations $A, B, Q, R$.

We will take the limit $t\to\infty$ in both sides of the identity
\eqref{eq:piid} with the scaling
\begin{equation}\label{eq:5.5}
  n = [2t + xt^{1/3}], \qquad z= -1 + \frac{w}{t^{1/3}}
\end{equation}
for a fixed real number $x$ and a complex number $w$, where $[a]$
denotes the largest integer smaller than or equal to $a$. We will
see that under this scaling limit, the right-hand-side of the
identity \eqref{eq:piid} becomes $f(x,w)$ given in \eqref{eq:fdef}
and the left-hand-side becomes~\eqref{eq:fPain}, thereby yielding
the desired Painlev\'e formula of $f(x,w)$.

Indeed, the limit of $e^{tz}\pi_n^*(z)$ is obtained in \cite{BR2}
and also \cite{BR4}. The paper \cite{BR2} studied the asymptotic
behavior of the longest increasing subsequences of certain
symmetrized versions of permutations and the asymptotic analysis
of $e^{tz}\pi_n^*(z)$ was a technical part of the paper. The paper
\cite{BR4} on the other hand studied a last passage percolation
model, which is different from the one discussed in Section
~\ref{sec:models}. From (5.22) and (5.26) of Proposition 5.4 and
Corollary 5.5 of \cite{BR2}, we find that
\begin{equation}
  e^{tz}\pi^*_n(z) \to
  \begin{cases} M_{22}(-\frac12iw;x), \qquad & w>0 \\
  -M_{21}(-\frac12iw;x)e^{\frac13 w^3-xw}, \qquad & w<0.
  \end{cases}
\end{equation}
This result actually motivated us to use the
function~\eqref{eq:symbolsp}.

On the other hand, it is known that \cite{BOO, kurtj:disc} (see
also \cite{WidomMoments})
\begin{equation}
  P_nABP_n \to \mathbf{A}_x
\end{equation}
in trace norm for any fixed real number $x$ where $\mathbf{A}_x$
is the Airy operator defined in \eqref{eq:Airyoperator}. This
limit was studied in the papers \cite{BOO, kurtj:disc,
WidomMoments} in the context of the longest increasing
subsequences and the Plancherel measure on partitions. Therefore,
the only remaining part is the asymptotic analysis of $Q$ and $R$.
These can be done by a standard steepest-descent analysis. Similar
analysis appeared in several places (see e.g. \cite{BOO, GTW,
BBP}) and we only sketch basic ideas.

We only consider $R$ since the analysis of $Q$ is similar. Note
that the integral formula \eqref{eq:defQR}, which was originally
defined for $|z|<1$, can be analytically continued for all complex
numbers $z$ by deforming the contour so that $z$ lies inside the
contour. To compute the limit of the right-hand-side of
\eqref{eq:piid}, we need the limit of $R([2t+yt^{1/3}])$ with
certain uniformity for $y\in [x,\infty)$ to ensure the convergence
of the inner product. As one can check from the analysis, it is
reasonable to think that $R([2t+yt^{1/3}])$ is close to
$R(2t+yt^{1/3})$ and we will compute the limit of the later. See
\cite{GTW}, for example, for a discussion of this type. Now from
\eqref{eq:5.3} and \eqref{eq:5.5}
\begin{equation}\label{eq:asproof1}
  R(2t+yt^{1/3}) = \frac{1}{2\pi i} \int
  \frac{-1+wt^{-1/3}}{b-(-1+wt^{-1/3})} (-b)^{2t+yt^{1/3}} e^{t(b-b^{-1})} da
\end{equation}
where the contour is modified to go from $\infty+i0$ to
$\infty-i0$ enclosing the origin and the point $-1+wt^{-1/3}$.
Here the function $(-b)^{2t+yt^{1/3}}$ denotes the principal
branch. Note that the integrand is of the form
\begin{equation}
  \frac{-1+wt^{-1/3}}{b-(-1+wt^{-1/3})} (-b)^{yt^{1/3}} e^{tf(b)}
\end{equation}
where
\begin{equation}
  f(b) = 2\log(-b)+b-b^{-1}.
\end{equation}
The function $f(b)$ has the double critical point at $b=-1$, and
$f(-1)=f'(-1)=f''(-1)=0$ and $f^{(3)}(-1)=2$. Approximately the
steepest-descent contour passing the critical point $b=-1$ is, in
a neighborhood of size, say $\epsilon>0$, of $b=-1$, the union of
the line from $-1+\epsilon e^{\pi i/3}$ to $-1$ and its complex
conjugate. As the pole $-1+w^{1/3}$ lies to the right of $b=-1$
when $w>0$, one can also check that it is possible to deform the
original contour to the steepest-descent contour when $w>0$. When
$w<0$, we can modify the contour to the union of $\{ -1+xe^{\pi
i/3} : 2|w|t^{-1/3}\le x\le \epsilon\}$,
$\{2|w|t^{-1/3}e^{i\theta} : \pi/3 \le\theta\le \pi\}$ and their
complex conjugate so that the pole $-1+wt^{-1/3}$ still lies on
the right of the contour but the contour `essentially' passes the
point $b=-1$. See \cite{BBP} where a similar modification of the
contour was used for a steepest-descent analysis. Note that in the
neighborhood of $b=-1$, the contour is oriented from $-1+\epsilon
e^{\pi i/3}$ to $-1+\epsilon e^{-\pi i/3}$.

 From the standard steepest-descent method, the integral is
asymptotic to the integral over the part of the contour in the
$\epsilon$ neighborhood of $a=-1$. The approximation $f(b)\simeq
\frac1{3!} f^{(3)}(-1)(b+1)^3= -\frac13 (b+1)^3$ suggests the
change of variables $it^{1/3}(b+1)=s$, which implies that
\begin{equation}
\begin{split}
  R(2t+yt^{1/3}) & \simeq \frac1{2\pi i} \int
  \frac{-1+wt^{-1/3}}{(-1-ist^{-1/3}-(-1+wt^{-1/3})}
  (1+ist^{-1/3})^{yt^{1/3}} e^{i\frac13 s^3} (-it^{-1/3})ds \\
  & \simeq \frac{-1}{2\pi} \int \frac1{is+w} e^{iys+i\frac13 s^3}
  ds
\end{split}
\end{equation}
where the contour is from $\infty e^{5\pi i/6}$ to $\infty e^{\pi
i/6}$ such that the pole $s=iw$ is above the contour. Hence we
find that $R(2t+yt^{1/3}) \simeq -C_w(y)$. Similar calculation
shows that $Q(2t+yt^{1/3}) \simeq -\Ai(y)$. This argument can be
made rigorous with uniform error bound for $y$ (see e.g.
\cite{BBP} for a similar calculation). Therefore, by noting that
$1-\mathbf{A}_x$ is a self-adjoint operator, one finds that the
right-hand-side of \eqref{eq:piid} converges to \eqref{eq:fdef}.
Thus we obtain the identity
\begin{equation}
    1-<\frac1{1-\mathbf{A}_x} C_w, \Ai>_{L^2((x,\infty))}
  = \begin{cases} M_{22}(-\frac12iw;x), \qquad & w>0 \\
  -M_{21}(-\frac12iw;x)e^{\frac13 w^3-xw}, \qquad & w<0.
  \end{cases}
\end{equation}
The proof of Theorem~\ref{thm:fPain} is complete.


\end{document}